\documentclass[11pt]{article}

\usepackage{amscd,amsmath,amssymb}
\usepackage[mathscr]{eucal}
\usepackage[all]{xy}
\usepackage{mathpple}

\makeatletter \@addtoreset{equation}{section}\makeatother

\newtheorem{theorem}{Theorem}[section]

\newtheorem{proposition}[theorem]{Proposition}

\newtheorem{definition}[theorem]{Definition}

\newcommand{\Z}{{\mathbb{Z}}}

\newcommand{\Ch}{{\mathcal{HH}}_0}

\newcommand{\Hom}{{\mathrm{Hom}}}
\newcommand{\Id}{{\mathrm{Id}}}
\newcommand{\Tr}{{\mathrm{Tr}}}

\newcommand{\cC}{{\mathcal{C}}}

\newcommand{\cL}{{\mathcal{L}}}

\newcommand{\cT}{{\mathcal{T}}}

\title{\bf 2-Representations and Equivariant 2D Topological Field Theories}
\author{D. Shklyarov}

\begin{document}

\maketitle

\tableofcontents

\bigskip

\section{Introduction}
This paper grew out of the author's (continued) attempt to
understand extended topological field theories, in particular, the
so-called Baez-Dolan hypothesis, and possible generalizations of all
that to Turaev's equivariant topological field theories \cite{T}.

The idea of extending the classical definition of TFTs to include
higher-categorical phenomena goes back to the early 90's and was
mostly motivated by Chern-Simons theory \cite{F}. The subject
has been developed significantly since then. We won't attempt to
give any rigorous definitions in this Introduction. Instead, we will
explain what the content of this paper has to do with extended TFTs
and how the aforementioned Baez-Dolan hypothesis motivates our
results. The reader, familiar with the state of affairs in this
subject, will certainly observe that we present the story in an
oversimplified form. We refer the reader to \cite{BD,F1}
for more thorough treatments.

As its title indicates, this paper deals with 2-dimensional TFTs only.
So from now on, a ``TFT'' will mean a ``2-dimensional TFT''. Let us also
fix a ground field, $K$, once and for all.

According to M. Atiyah, a TFT is a rule that assigns a
finite-dimensional vector space $C$ to an oriented circle $S^1$, the
tensor power $C^{\otimes n}$ to the disjoint union of $n$ circles,
and linear maps between the tensor powers to (isomorphism classes
of) oriented 2-dimensional cobordisms between the unions of circles.
The axioms that this assignment is required to satisfy are most
conveniently expressed by saying that it is a monoidal functor from
the symmetric monoidal category ${\bf Cob_2}$, whose objects are
closed oriented 1-manifolds and morphisms are 2-cobordisms, to the
category of finite-dimensional vector spaces.

TFTs are well studied objects. It is a classical result that the
vector space $C$ a TFT assigns to the circle carries a {\it
commutative Frobenius algebra} structure, and conversely, any such
algebra determines a TFT.

The definition of TFT can be refined as follows \cite{BD}: An
extended TFT is a rule that assigns small linear categories to
oriented 0-manifolds, linear functors to oriented 1-cobordisms of
0-manifolds, and natural transformations to oriented 2-cobordisms
with corners. As before, the assignment should give rise to a
monoidal 2-functor from the symmetric monoidal 2-category ${\bf
Cob_2^{ext}}$, whose objects are 0-manifolds, 1-morphisms are
1-cobordisms, and 2-morphisms are 2-cobordisms between 1-cobordisms,
to the 2-category of small linear categories. Notice that an
extended TFT gives rise to a TFT in the usual sense.

The above definition can be spelled out for other ``target''
2-categories. We will be interested in the 2-category ${\bf Bim}_K$
whose objects are algebras, 1-morphisms are bimodules, and
2-morphisms are morphisms of bimodules. This 2-category is a first
approximation to what should be called ``the 2-category of
non-commutative affine schemes over $K$''. This non-commutative
geometric interpretation is quite relevant in this setting as we
will see later.

There is, in fact, an analog for extended TFTs of the explicit
description of the usual TFTs in terms of commutative Frobenius
algebras. Namely, J. Baez and J. Dolan have conjectured \cite{BD}
that extended TFTs, valued in a $2$-category ${\bf C}$, should be
completely determined by the objects they assign to an oriented
point; moreover, they have observed that the image of the point
should be dualizable in ${\bf C}$ in a certain strong sense, and
vice versa, every dualizable object in ${\bf C}$ should give rise to
an extended TFT.

A result of this sort has been proved recently by Hopkins and Lurie
(see e.g. \cite[Section 3]{F1}). The result implies that in the case
${\bf C}={\bf Bim}_K$ the extended TFTs are, roughly speaking, in
one-to-one correspondence with {\it symmetric Frobenius separable
algebras}, i.e. separable algebras equipped with a non-degenerate
trace\footnote{From now on, by a trace we will understand a
cyclically invariant functional, so we will omit the word
``symmetric''.}.

This result admits the following non-commutative geometric
interpretation: essentially, it suggests that extended TFTs valued
in the ``right'' 2-category of non-commutative schemes should be in
one-to-one correspondence with non-commutative smooth compact
Calabi-Yau spaces. Apparently, this more general statement also
follows from the results of Hopkins and Lurie. This is in perfect
agreement with the result of K. Costello \cite{Cos} establishing
equivalence between a category of non-commutative Calabi-Yau spaces
and a category of open-closed chain-level TFTs (see also \cite{KS}).

We believe that the idea of extended TFTs, as well as the language
of non-commutative geometry, will prove useful beyond the above
setting. For instance, one can try to apply these ideas to the
so-called homotopy field theories \cite{BT,BTW,BuTW,R,PT,T} in which
manifolds and cobordisms are decorated by maps to a target space.

One of the first cases to look at is the case when the target is the
classifying space of a finite group, $G$ \cite{T}. In this case
homotopy field theories (a.k.a. $G$-equivariant TFTs \cite{MS}) are
symmetric monoidal functors from the category of principal
$G$-bundles over manifolds, with morphisms being principal
$G$-bundles over cobordisms, to the category of vector spaces.
Equivariant TFTs provide an adequate language for describing the
orbifolding procedure in the setting of TFTs with symmetries
\cite{K1,K2,MS}.

We will not attempt to define {\it extended} equivariant TFTs in
this paper. Instead, we would like to emphasize the non-commutative
geometric aspect of the sought-for theory. Namely, we believe that
the right definition, whatever it is, should lead to the following
statement: extended $G$-equivariant TFTs valued in the 2-category of
non-commutative schemes are in one-to-one correspondence with
non-commutative smooth compact Calabi-Yau spaces acted on by $G$. In
particular, extended $G$-equivariant TFTs valued in ${\bf C}={\bf
Bim}_K$ should be in one-to-one correspondence with {\it Frobenius
separable algebras equipped with a categorical $G$-action preserving
the Frobenius structure}. The meaning of the adjective
``categorical'' will be explained in the next section.

The aim of this paper is to present a piece of evidence in favor of
the above point of view:
\begin{quote}
{\bf Main result:} {\it We explicitly construct a (non-extended)
$G$-equivariant TFT from an arbitrary, not necessarily separable
Frobenius algebra with a categorical $G$-action. In such a TFT, the
connected components of cobordisms are required to have at least one
incoming boundary component (cf. \cite{Cos, KS}). When the algebra
is separable, we get an honest equivariant TFT.}
\end{quote}

We note that the idea of producing equivariant TFTs from Frobenius
algebras equipped with a $G$-action is not new: it arises naturally
in the study of open-closed equivariant TFTs \cite[Section 7]{MS}.
The construction offered in the present paper is, in a sense, dual
to the one presented in \cite{MS} (the meaning of ``dual'' will be
explained in Section \ref{conclrem}). Another difference is that we
work in the more general setting of non-separable Frobenius algebras
and {\it twisted}(=categorical) $G$-actions on them.

For the reader's convenience, we introduce all the necessary
definitions and formulate the results in a concise way in the first
part of the paper. The proofs are collected in appendices in the
remaining part. Throughout this paper, $G$ and $K$ stand for a
finite group and a field, respectively.

\medskip

{\bf Acknowledgements.} I am grateful to K. Costello, D. Freed, and
Y. Soibelman for various useful remarks on the results and the
exposition. I would also like to thank the authors of \cite{GK,TV,W}
and especially \cite{MS} for inspiration.

\medskip

\section{Frobenius algebras with twisted $G$-action}\label{cr}
Let us begin by recalling the notion of categorical representation
of a group \cite{GK}.

A categorical representation of a group $G$ is a category $\cC$
together with the following data:
  \begin{enumerate}
    \item for each element $g\in G$, an autoequivalence $\rho(g)$ of
      $\cC$;
    \item for any pair of elements $g,h\in G$, an isomorphism of functors
      $$
        c(g,h): \rho(g)\circ\rho(h)\Longrightarrow\rho(gh);
      $$
    \item an isomorphism of functors
      $$
        c(e): \rho(e)\Longrightarrow \Id_{\cC},
      $$
      where $e$ is the unit element of $G$.
  \end{enumerate}
  The above autoequivalences and isomorphisms are required to satisfy
  the following conditions:
  \begin{enumerate}
      \item for any triple of elements $g,h,k\in G$, the diagram
      \begin{displaymath}
      \xymatrix{ \rho(g)\circ\rho(h)\circ\rho(k)
      \ar@{=>}[r]^{\qquad\rho(g)\circ
      c(h,k)}\ar@{=>}[d]_{c(g,h)\circ\rho(k)} &
      \quad\rho(g)\circ\rho(hk) \ar@{=>}[d]^{c(g,hk)}\\
      \rho(gh)\circ\rho(k) \ar@{=>}[r]^{c(gh,k)} & \rho(ghk)}
      \end{displaymath}
      commutes;
    \item for any element $g\in G$, the diagrams
    \begin{displaymath}
      \xymatrix{ \rho(g)\circ\rho(e) \ar@{=>}[r]^{\quad\rho(g)\circ
      c(e)}\ar@{=>}[rd]_{c(g,e)} &
      \,\,\rho(g)\circ\Id_{\cC} \ar@{=}[d]\\
       & \rho(g)} \qquad \xymatrix{ \rho(e)\circ\rho(g)
      \ar@{=>}[r]^{\,\,\, c(e)\circ\rho(g)}\ar@{=>}[rd]_{c(e,g)}
      &
      \,\,\Id_{\cC}\circ\rho(g) \ar@{=}[d]\\
       &
      \rho(g)}
      \end{displaymath}
      commute.
     \end{enumerate}

From now on, we will only consider the case when the category $\cC$ is
$K$-linear and has one object. We will identify such a category
with its endomorphism algebra. In this case, $\rho(g)$ are
automorphisms of $\cC$, whereas the isomorphisms $c(g,h)$ and $c(e)$
are given by conjugation with some invertible elements of $\cC$, which we denote by
$c_{g,h}$ and $c_e$. Namely, the latter are defined by
      $$
        \rho(g)\cdot\rho(h)=Ad(c_{g,h})\cdot\rho(gh),\quad \rho(e)=Ad(c_e),
      $$
where $Ad(c)$ stands for the inner automorphism $c_1\mapsto c\cdot
c_1\cdot c^{-1}$ and the equalities are understood as equalities of
automorphisms of $\cC$. The above commutative diagrams boil down to
the following equalities
\begin{equation}\label{cocycle}
c_{g,h}c_{gh,k}=g(c_{h,k})c_{g,hk},\quad c_{g,e}=g(c_e),\quad c_{e,g}=c_e,
\end{equation}
where we write $g(c)$ instead of $\rho(g)(c)$.

\medskip

\begin{definition}
A {\bf twisted algebra bundle} on the classifying stack $BG$ is an
algebra $\cC$ equipped with a categorical $G$-action as above.
\end{definition}

\medskip

Later on, we will need one more definition. Recall that a Frobenius
algebra is a finite-dimensional (not necessarily commutative) unital
algebra $\cC$ equipped with a trace $\theta: \cC\to K$ such that the
pairing $(c_1,c_2)\mapsto \theta(c_1c_2)$ is non-degenerate.

\medskip

\begin{definition}
A {\bf twisted Frobenius algebra bundle} on the classifying stack
$BG$ is a Frobenius algebra $(\cC,\theta)$ equipped with a categorical
$G$-action satisfying $\theta(g(c))=\theta(c)$ for any $g\in G$ and
$c\in\cC$.
\end{definition}

\medskip

\section{Equivariant topological field theories}\label{etft}
Turaev has shown \cite{T} that the data of a $G$-equivariant TFT is
equivalent to that of a crossed $G$-algebra, a $G$-equivariant
analog of a commutative Frobenius algebra. In this section, we
recall Turaev's definition of crossed $G$-algebra. It will be
convenient for us to first introduce some auxiliary notions.

\medskip

\begin{definition}\label{bun} A {\bf vector bundle on the loop space} $\cL BG$ is a $G$-graded
finite-dimensional vector space $C=\oplus_g C_g$ together with a
group homomorphism $G\to GL(C)$ such that $g(C_h)=C_{ghg^{-1}}$
for all $g,h\in G$.
\end{definition}
\begin{definition}
A {\bf special\footnote{Apparently, such objects should be called ``$S^1$-equivariant bundles on $\cL BG$''. We decided to stick to the neutral term ``special'' to avoid the necessity of justifying our terminology.} vector bundle on $\cL BG$} is a vector bundle such that $g|_{C_g}={\rm id}$.
\end{definition}

\medskip

Now we are in position to define crossed $G$-algebras.

\begin{definition}\label{def1} A {\bf crossed
$G$-algebra} (or {\bf Turaev algebra} in the terminology of \cite{MS}) is a
special vector bundle $C=\oplus_g C_g$ on $\cL BG$ together with a $G$-invariant functional
$\theta_e:C_e \to K$ and an algebra structure satisfying the following properties:
\begin{itemize}
\item[(1)] $G$ acts by automorphisms of the algebra $C$;
\item[(2)] $C_g\cdot C_h\subset C_{gh}$;
\item[(3)] for all $c'\in C_g$ and $c''\in C_h$, $$c'c''=g(c'')c';$$
\item[(4)] (Torus axiom) for all $g,h\in G$ and $c\in C_{hgh^{-1}g^{-1}}$,
  $$
  {\Tr}_{C_h}(L'_c\cdot g)= {\Tr}_{C_{g}}(h^{-1}\cdot L''_c),
  $$
where $L'_c:C_{ghg^{-1}}\to C_{h}$ and $L''_c: C_{g}\to
C_{hgh^{-1}}$ stand for the operators of left multiplication
with $c$.
\item[(5*)] $C$ is unital;
\item[(6*)] $\theta_e$ induces a non-degenerate pairing $C_g \otimes C_{g^{-1}}\to K$.
\end{itemize}
\end{definition}

\medskip

We will also need a weaker notion which we call a weak crossed
$G$-algebra. Weak crossed $G$-algebras correspond to equivariant
TFTs in which the connected components of cobordisms are required to
have at least one incoming boundary component.

\begin{definition}\label{def2}
A {\bf weak crossed $G$-algebra} is an object $C$ satisfying all the
requirements listed in Definition \ref{def1} {\bf except for} (5*)
and (6*); instead, $C$ possesses the following extra structure:
there is a coassociative coalgebra structure $\Delta: C\to C\otimes
C$ satisfying the following properties
\begin{itemize}
\item[(5)] $G$ acts by automorphisms of the coalgebra $C$;
\item[(6)] $\Delta$ respects the $G$-grading, i.e. $\Delta(C_k)\subset\oplus_{gh=k}\,C_{g}\otimes C_{h}$ (we will denote the corresponding map $C_{gh}\to C_{g}\otimes C_{h}$ by $\Delta_{g,h}$);
\item[(7)] for all $g,h\in G$, $\Delta_{g,h}=\sigma (1\otimes h)\Delta_{h,h^{-1}gh}$ (here $\sigma$ is the transposition map);
\item[(8)] for any $g\in G$, $(\theta_e\otimes1)\Delta_{e,g}=(1\otimes\theta_e)\Delta_{g,e}={\rm id}_{C_g}$;
\item[(9)] $\Delta$ is a morphism of $C$-bimodules.
\end{itemize}
\end{definition}

Note that any crossed $G$-algebra possesses a weak crossed
$G$-algebra structure: the maps $\Delta_{g,h}$ come from the
morphisms in the corresponding equivariant TFT defined by principal
$G$-bundles on the genus 0 surface with one incoming and two
outgoing boundaries.

Also observe that a weak crossed $G$-algebra is an honest crossed
$G$-algebra iff it is unital. Indeed, if $1_C$ is the unit then for
any $c\in C_g$ we have
$$
(1\otimes\theta_e)(1\otimes m_{g^{-1},g})(\Delta_{g,g^{-1}}(1_C)\otimes c)\stackrel{(9)}=(1\otimes\theta_e)(\Delta_{g,1}(c))\stackrel{(8)}=c
$$
which implies triviality of the kernel of the pairing, defined by $\theta_e$.

\medskip

\section{A $G$-equivariant version of the 0-th Hochschild homology}\label{cHochschild}

In this section, we will introduce a $G$-equivariant version of the 0-th Hochschild homology
of an algebra.

Consider a twisted algebra bundle, i.e. an algebra $\cC$ equipped
with a categorical $G$-action (we will keep the notations from
Section \ref{cr}). For an element $h\in G$ set
$$
\cC_h=\text{span}\{c_1c_2-c_2h(c_1)\,|\,c_1,c_2\in \cC\}\subset\cC
$$
and define
$$
\Ch(\cC)=\oplus_g \Ch(\cC)_g, \quad \Ch(\cC)_g:=\cC/\cC_g.
$$

The space $\Ch(\cC)$ inherits a $G$-action from $\cC$. In order to describe it, we
need some auxiliary definitions and results.

For a pair of elements $g,h\in G$, define a linear map $T_h(g):\cC\to\cC$ by the
formula
$$
T_h(g): c\mapsto c_e^{-1}c_{g,g^{-1}}^{-1}g(c)c_{g,h}c_{gh,g^{-1}}.
$$

\medskip

\begin{proposition}\label{maps} For any $g,h\in G$, $T_h(g)(\cC_h)=\cC_{ghg^{-1}}$.
\end{proposition}
\noindent {\bf Proof} can be found in Appendix \ref{a1}.

\medskip

By the above proposition, we have linear maps
$$
\cT_h(g): \Ch(\cC)_h\to\Ch(\cC)_{ghg^{-1}}
$$
induced by the maps $T_h(g)$. Define $\cT(g): \Ch(\cC)\to\Ch(\cC)$ by
\begin{equation}\label{G}
\cT(g)|_{\Ch(\cC)_h}=\cT_h(g).
\end{equation}

\medskip

\begin{proposition}\label{repr}
The operators $\cT(g)$ form a representation of $G$ in $\Ch(\cC)$.
Moreover, $$\cT(g)|_{\Ch(\cC)_g}={\rm id}.$$
\end{proposition}
\noindent {\bf Proof} can be found in Appendix \ref{a2}.

\medskip

Now we are ready to formulate
\begin{definition}
The special vector bundle $\Ch(\cC)$ on $\cL BG$ will be called {\bf
the 0-th Hochschild homology bundle} of $\cC$.
\end{definition}


\medskip

\section{An equivariant TFT structure on the 0-th Hochschild homology bundle}\label{etftfromcr}

In this section, we equip the 0-th Hochschild homology bundle of a twisted
Frobenius algebra bundle with a weak crossed $G$-algebra structure.

Let us fix a twisted Frobenius algebra bundle, i.e. a Frobenius
algebra $\cC=(\cC,\theta)$ equipped with a categorical $G$-action
such that $\theta$ is $G$-invariant (see Section \ref{cr}). Let
$\xi=\sum_i\xi'_i\otimes \xi''_i\in \cC\otimes\cC$ stand for the
symmetric tensor inverse to the pairing defined by $\theta$:
$$
c=\sum_i \xi'_i\theta(c\xi''_i),\quad \forall c\in\cC.
$$
(In what follows, we omit the summation and write $\xi'_i\otimes \xi''_i$.)

For any pair of elements $g,h\in G$ define linear maps $m_{g,h}:\cC\otimes\cC\to\cC$ by
\begin{equation}\label{m}
m_{g,h}(c'\otimes c'')=\xi_i'c'g(\xi_i''c'')c_{g,h}.
\end{equation}
\begin{proposition}\label{alg}
The maps $m_{g,h}$ descend to well-defined linear maps $$m_{g,h}:
\Ch(\cC)_g\otimes\Ch(\cC)_h\to\Ch(\cC)_{gh}.$$ Altogether, they
define an associative algebra structure on $\Ch(\cC)$.
\end{proposition}
\noindent {\bf Proof} can be found in Appendix \ref{a3}.
\medskip

Now for any pair of elements $g,h\in G$ define linear maps
$\Delta_{g,h}:\cC\to\cC\otimes\cC$ by
\begin{equation}\label{com}
\Delta_{g,h}(c)=c\,c_{g,h}^{-1}g(\xi_i')\otimes \xi_i''.
\end{equation}
\begin{proposition}\label{coalg}
The maps $\Delta_{g,h}$ descend to well-defined linear maps
$$\Delta_{g,h}: \Ch(\cC)_{gh}\to\Ch(\cC)_g\otimes\Ch(\cC)_h.$$
Altogether, they define a coassociative coalgebra structure on
$\Ch(\cC)$.
\end{proposition}
\noindent {\bf Proof} can be found in Appendix \ref{a3}.

\medskip

Let us formulate the main result of this paper:
\begin{theorem}\label{main}
For any twisted Frobenius algebra bundle $\cC=(\cC,\theta)$ on $BG$,
the 0-th Hochschild homology bundle $\Ch(\cC)$, equipped with the $G$-action
(\ref{G}), the multiplication (\ref{m}), the comultiplication
(\ref{com}), and the functional
\begin{equation}\label{theta}
\theta_e: \Ch(\cC)_e\to K, \quad \theta_e(c)=\theta(cc_e)
\end{equation}
is a weak crossed $G$-algebra.
\end{theorem}
\noindent {\bf Proof} can be found in Appendix \ref{a4}.

\medskip

Using the observation, mentioned at the end of Section \ref{etft},
one can show that $\Ch(\cC)$ is an honest crossed $G$-algebra iff
the central element $\xi_i'\xi_i''$ is invertible\footnote{This is
so when $\cC$ is separable.}. Indeed, in this case
$(\xi_i'\xi_i''c_e)^{-1}$ is the unit of $\Ch(\cC)$, as one can see
from the following computation:
$$
m_{e,g}((\xi_i'\xi_i''c_e)^{-1}\otimes c)=\xi_j'(\xi_i'\xi_i''c_e)^{-1}e(\xi_j''c)c_{e,g}=\xi_j'(\xi_i'\xi_i'')^{-1}\xi_j''c=c,
$$
$$
m_{g,e}(c\otimes (\xi_i'\xi_i''c_e)^{-1})=\xi_j'cg(\xi_j''(\xi_i'\xi_i''c_e)^{-1})c_{g,e}=\xi_j'cg(\xi_j''(\xi_i'\xi_i'')^{-1})g(c_e)^{-1}c_{g,e}
$$
$$
=\xi_j'cg(\xi_j''(\xi_i'\xi_i'')^{-1})\stackrel{modulo\,\cC_g}=\xi_j''(\xi_i'\xi_i'')^{-1}\xi_j'c=c.
$$

\medskip

\section{Concluding remarks}\label{conclrem}

We already mentioned in the Introduction that the content of this
paper is closely related to that of Section 7 of \cite{MS}. In the
first part of this section we will expand a little bit on how our
result compares with those obtained in \cite{MS}. The second part is
devoted to a simple observation relating equivariant TFTs to
generalized group characters \cite{HKR1,HKR2, GK}.

In Section 7.3 of \cite{MS}, the authors describe a construction,
which goes back to \cite{T}, that associates a $G$-equivariant TFT
with a finite $G$-space, $X$, equipped with a $G$-invariant ``volume
form'' and a $B$-field. The latter is an element of the second
cohomology group $H^2(G,A(X)^\times)$ of $G$ with values in the
abelian group $A(X)^\times$ of invertible functions on $X$. Such
data -- a $G$-space with an invariant trace and a $B$-field -- give
rise to a Frobenius algebra bundle on $BG$ in the sense of the
present paper. Notice that the ``fiber'' of this bundle is a
commutative semisimple algebra. Reversing the logic, one can say
that our result is about non-commutative $G$-spaces\footnote{The
elements $c_{g,h}$ should then be viewed as defining a $B$-field
which belongs to some ``non-abelian cohomology group''
$H^2(G,\cC^\times)$.}.

Actually, a special case of non-commutative $G$-spaces is mentioned
already in \cite{MS}. Namely, in Section 7.4, where the authors
introduce and study open-closed equivariant TFTs, they describe an
explicit construction of equivariant TFTs from Frobenius semisimple
categories with $G$-action (see Theorem 10 in loc. cit.). However,
by a category with $G$-action the authors understand a category
enriched in the category of representations of $G$. It means, in
particular, that one has an {\it honest}, non-twisted action of $G$
on the morphism spaces. On a final note, our construction of
equivariant TFTs differs from the one described in \cite{MS} even in
the cases when both constructions are applicable. In our approach,
the underlying space of the equivariant TFT is a $G$-equivariant
version of the 0-th Hochschild {\it homology} of an algebra whereas
in \cite{MS} it is a $G$-equivariant version of the
center\footnote{Observe that this is precisely the categorical trace
in the sense of \cite{GK,B}.} (=the 0-th Hochschild {\it
cohomology}). The reason we have decided to work with the
equivariant Hochschild homology is that the latter, we believe, is
more suitable for the purpose of constructing equivariant extended
TFTs.

We would like to conclude this section with the following curious
(although straightforward) observation. In their recent work
\cite{GK}, N. Ganter and M. Kapranov have emphasized the importance
of categorical representations of groups in the study of generalized
equivariant cohomology theories and related subjects. In particular,
they have noticed that the so-called 2-class functions on finite
groups studied in \cite{HKR1,HKR2} arise naturally as the
2-characters of 2-representations of the groups. Namely, one defines
the categorical character of an arbitrary 2-representation which is
a vector bundle $\mathbf{Tr}=\oplus_g \mathbf{Tr}_g$ on $\cL BG$ in
the sense of Definition \ref{bun}; then the 2-character is the
``character'' of the categorical character:
$$
\chi(g,h)={\Tr}(g: \mathbf{Tr}_h\to \mathbf{Tr}_h)
$$
where $gh=hg$. Clearly, $\chi(kgk^{-1},khk^{-1})=\chi(g,h)$ for any
$k\in G$ which is what being a {\it 2-class function} means.

The point we would like to stress is that by the above mentioned
results of \cite{MS} the categorical traces of 2-representations of
a finite group $G$ in 2-vector spaces can be upgraded to
$G$-equivariant TFTs. This implies that the 2-characters are {\it
modular invariant} in the following sense: for any $\begin{pmatrix}
a & b\\c & d\end{pmatrix}\in SL(2,\Z)$ one has
$$
\chi(g^ah^b,g^ch^d)=\chi(g,h)
$$
for any pair of commuting elements $g,h\in G$. Indeed, it is enough
to show that $\chi$ is preserved by the generators $\begin{pmatrix}
1 & 1\\0 & 1\end{pmatrix}$ and $\begin{pmatrix} 0 & -1\\1 &
0\end{pmatrix}$ of $SL(2,\Z)$, i.e.
$$
\chi(gh,h)=\chi(g,h),\quad \chi(h^{-1},g)=\chi(g,h).
$$
The first equality follows from the fact that $\mathbf{Tr}$ is a
special vector bundle, i.e. $h|_{\mathbf{Tr}_h}={\rm id}$; the
second equality is a consequence of the torus axiom.

\pagebreak

\appendix
\section{Proof of Proposition \ref{maps}}\label{a1}
Let us compute $T_h(g)(c_1c_2-c_2h(c_1))$:
$$
T_h(g)(c_1c_2-c_2h(c_1))=c_e^{-1}c_{g,g^{-1}}^{-1}\left(g(c_1)g(c_2)-g(c_2)g(h(c_1))\right)c_{g,h}c_{gh,g^{-1}}
$$
$$
=c_e^{-1}c_{g,g^{-1}}^{-1}\left(g(c_1)g(c_2)- g(c_2)Ad(c_{g,h})(gh(c_1))\right)c_{g,h}c_{gh,g^{-1}}
$$
$$
=c_e^{-1}c_{g,g^{-1}}^{-1}\left(g(c_1)g(c_2)- g(c_2) Ad(c_{g,h})Ad(c_{ghg^{-1},g}^{-1})(ghg^{-1}(g(c_1)))\right)c_{g,h}c_{gh,g^{-1}}
$$
$$
=c_e^{-1}c_{g,g^{-1}}^{-1}\left(c'c''- c''Ad(c_{g,h})Ad(c_{ghg^{-1},g}^{-1})(ghg^{-1}(c'))\right)c_{g,h}c_{gh,g^{-1}}
$$
\begin{equation}\label{1}
=c_e^{-1}c_{g,g^{-1}}^{-1}c'c''c_{g,h}c_{gh,g^{-1}}-c_e^{-1}c_{g,g^{-1}}^{-1}c'' c_{g,h}c_{ghg^{-1},g}^{-1}ghg^{-1}(c')c_{ghg^{-1},g}c_{gh,g^{-1}},
\end{equation}
where $c'=g(c_1)$, $c''=g(c_2)$. Observe that, modulo $\cC_{ghg^{-1}}$, the first summand in (\ref{1})
is equal to
$$
c'c''c_{g,h}c_{gh,g^{-1}}ghg^{-1}(c_e^{-1})ghg^{-1}(c_{g,g^{-1}}^{-1}),
$$
which, in its turn, is equal to
$$
c'c''c_{g,h}c_{ghg^{-1},g}^{-1}
$$
due to (\ref{cocycle}). On the other hand, modulo $\cC_{ghg^{-1}}$, the second summand in (\ref{1}) is equal to
$$
c'' c_{g,h}c_{ghg^{-1},g}^{-1}ghg^{-1}(c')c_{ghg^{-1},g}c_{gh,g^{-1}}ghg^{-1}(c_e^{-1})ghg^{-1}(c_{g,g^{-1}}^{-1})
$$
and, by (\ref{cocycle}), this is the same thing as
$$
c'' c_{g,h}c_{ghg^{-1},g}^{-1}ghg^{-1}(c').
$$
Thus, modulo $\cC_{ghg^{-1}}$, the expression (\ref{1}) is equal to
$$
c'(c''c_{g,h}c_{ghg^{-1},g}^{-1})-(c'' c_{g,h}c_{ghg^{-1},g}^{-1})ghg^{-1}(c')\equiv0.
$$
The proposition is proved.


\section{Proof of Proposition \ref{repr}}\label{a2}
Let us apply $\cT(g_1)\cT(g_2)$ and $\cT(g_1g_2)$ to an element $c\in\Ch(\cC)_h$:
$$
\cT(g_1)\cT(g_2)(c)=\cT(g_1)(c_e^{-1}c_{g_2,g_2^{-1}}^{-1}g_2(c)c_{g_2,h}c_{g_2h,g_2^{-1}})
$$
$$
=c_e^{-1}c_{g_1,g_1^{-1}}^{-1}g_1(c_e^{-1}c_{g_2,g_2^{-1}}^{-1}g_2(c)c_{g_2,h}c_{g_2h,g_2^{-1}})c_{g_1,g_2hg_2^{-1}}c_{g_1g_2hg_2^{-1},g_1^{-1}}
$$
\begin{equation}\label{2}
=c_e^{-1}c_{g_1,g_1^{-1}}^{-1}g_1(c_e^{-1})g_1(c_{g_2,g_2^{-1}}^{-1})c_{g_1,g_2}g_1g_2(c)c_{g_1,g_2}^{-1}g_1(c_{g_2,h})g_1(c_{g_2h,g_2^{-1}})c_{g_1,g_2hg_2^{-1}}c_{g_1g_2hg_2^{-1},g_1^{-1}}.
\end{equation}
On the other hand,
\begin{equation}\label{3}
\cT(g_1g_2)(c)=c_e^{-1}c_{g_1g_2,g_2^{-1}g_1^{-1}}^{-1}g_1g_2(c)c_{g_1g_2,h}c_{g_1g_2h,g_2^{-1}g_1^{-1}}.
\end{equation}
We have to show that (\ref{2}) coincides with (\ref{3}) modulo $\cC_{g_1g_2hg_2^{-1}g_1^{-1}}$.
We will use the relation $c_1c_2-c_2g_1g_2hg_2^{-1}g_1^{-1}(c_1)$ to ``move'' everything that is on the left of $g_1g_2(c)$ in (\ref{2}) and (\ref{3}) to the right side of the expressions.

Thus, modulo $\cC_{g_1g_2hg_2^{-1}g_1^{-1}}$, the expressions (\ref{2}) and (\ref{3}) are equal to
$$
g_1g_2(c)c_{g_1,g_2}^{-1}g_1(c_{g_2,h})g_1(c_{g_2h,g_2^{-1}})c_{g_1,g_2hg_2^{-1}}c_{g_1g_2hg_2^{-1},g_1^{-1}}\,g_1g_2hg_2^{-1}g_1^{-1}(c_e^{-1}c_{g_1,g_1^{-1}}^{-1}g_1(c_e^{-1})g_1(c_{g_2,g_2^{-1}}^{-1})c_{g_1,g_2})
$$
and
$$
g_1g_2(c)c_{g_1g_2,h}c_{g_1g_2h,g_2^{-1}g_1^{-1}}\,g_1g_2hg_2^{-1}g_1^{-1}(c_e^{-1}c_{g_1g_2,g_2^{-1}g_1^{-1}}^{-1}),
$$
respectively. Therefore, it is enough to show that
$$
c_{g_1,g_2}^{-1}g_1(c_{g_2,h})g_1(c_{g_2h,g_2^{-1}})c_{g_1,g_2hg_2^{-1}}c_{g_1g_2hg_2^{-1},g_1^{-1}}\,g_1g_2hg_2^{-1}g_1^{-1}(c_e^{-1}c_{g_1,g_1^{-1}}^{-1}g_1(c_e^{-1})g_1(c_{g_2,g_2^{-1}}^{-1})c_{g_1,g_2})
$$
$$
=c_{g_1g_2,h}c_{g_1g_2h,g_2^{-1}g_1^{-1}}\,g_1g_2hg_2^{-1}g_1^{-1}(c_e^{-1}c_{g_1g_2,g_2^{-1}g_1^{-1}}^{-1}).
$$
Observe that by (\ref{cocycle}), the left hand side of the latter
equality can be simplified as follows (we will underline the places
in the formulas that are about to be changed):
$$
c_{g_1,g_2}^{-1}g_1(c_{g_2,h})g_1(c_{g_2h,g_2^{-1}})c_{g_1,g_2hg_2^{-1}}c_{g_1g_2hg_2^{-1},g_1^{-1}}\,g_1g_2hg_2^{-1}g_1^{-1}(c_e^{-1}c_{g_1,g_1^{-1}}^{-1}\underline{g_1(c_e^{-1})g_1(c_{g_2,g_2^{-1}}^{-1})c_{g_1,g_2}})
$$
$$
=c_{g_1,g_2}^{-1}g_1(c_{g_2,h})g_1(c_{g_2h,g_2^{-1}})c_{g_1,g_2hg_2^{-1}}c_{g_1g_2hg_2^{-1},g_1^{-1}}\underline{\,g_1g_2hg_2^{-1}g_1^{-1}(c_e^{-1}c_{g_1,g_1^{-1}}^{-1}}c_{g_1g_2,g_2^{-1}}^{-1})
$$
$$
=c_{g_1,g_2}^{-1}g_1(c_{g_2,h})g_1(c_{g_2h,g_2^{-1}})c_{g_1,g_2hg_2^{-1}}c_{g_1g_2hg_2^{-1}g_1^{-1},g_1}^{-1}\,g_1g_2hg_2^{-1}g_1^{-1}(c_{g_1g_2,g_2^{-1}}^{-1}),
$$
so now we have to prove that
$$
c_{g_1,g_2}^{-1}g_1(c_{g_2,h})g_1(c_{g_2h,g_2^{-1}})c_{g_1,g_2hg_2^{-1}}c_{g_1g_2hg_2^{-1}g_1^{-1},g_1}^{-1}\,g_1g_2hg_2^{-1}g_1^{-1}(c_{g_1g_2,g_2^{-1}}^{-1})
$$
$$
=c_{g_1g_2,h}c_{g_1g_2h,g_2^{-1}g_1^{-1}}\,g_1g_2hg_2^{-1}g_1^{-1}(c_e^{-1}c_{g_1g_2,g_2^{-1}g_1^{-1}}^{-1}).
$$
By (\ref{cocycle}), the right hand side of the latter equality can be simplified as follows
$$
c_{g_1g_2,h}c_{g_1g_2h,g_2^{-1}g_1^{-1}}\underline{\,g_1g_2hg_2^{-1}g_1^{-1}(c_e^{-1}c_{g_1g_2,g_2^{-1}g_1^{-1}}^{-1})}
$$
$$
=c_{g_1g_2,h}c_{g_1g_2h,g_2^{-1}g_1^{-1}}c_{g_1g_2h,g_2^{-1}g_1^{-1}}^{-1}c_{g_1g_2hg_2^{-1}g_1^{-1},g_1g_2}^{-1}
$$
and now we will be proving that
$$
c_{g_1,g_2}^{-1}g_1(c_{g_2,h})g_1(c_{g_2h,g_2^{-1}})c_{g_1,g_2hg_2^{-1}}c_{g_1g_2hg_2^{-1}g_1^{-1},g_1}^{-1}\,g_1g_2hg_2^{-1}g_1^{-1}(c_{g_1g_2,g_2^{-1}}^{-1})
$$
\begin{equation}\label{5}
=c_{g_1g_2,h}c_{g_1g_2hg_2^{-1}g_1^{-1},g_1g_2}^{-1}.
\end{equation}
Let us simplify the left-hand side further:
$$
c_{g_1,g_2}^{-1}g_1(c_{g_2,h})g_1(c_{g_2h,g_2^{-1}})c_{g_1,g_2hg_2^{-1}}\underline{c_{g_1g_2hg_2^{-1}g_1^{-1},g_1}^{-1}\,g_1g_2hg_2^{-1}g_1^{-1}(c_{g_1g_2,g_2^{-1}}^{-1})}
$$
$$
\stackrel{by (\ref{cocycle})}=c_{g_1,g_2}^{-1}g_1(c_{g_2,h})g_1(c_{g_2h,g_2^{-1}})c_{g_1,g_2hg_2^{-1}}c_{g_1g_2h,g_2^{-1}}^{-1}c_{g_1g_2hg_2^{-1}g_1^{-1},g_1g_2}^{-1},
$$
so (\ref{5}) is equivalent to
$$
c_{g_1,g_2}^{-1}g_1(c_{g_2,h})g_1(c_{g_2h,g_2^{-1}})c_{g_1,g_2hg_2^{-1}}c_{g_1g_2h,g_2^{-1}}^{-1}=c_{g_1g_2,h}
$$
or
\begin{equation}\label{6}
c_{g_1,g_2}^{-1}g_1(c_{g_2,h}c_{g_2h,g_2^{-1}})c_{g_1,g_2hg_2^{-1}}=c_{g_1g_2,h}c_{g_1g_2h,g_2^{-1}}.
\end{equation}
Consider the left-hand side of the latter equality:
$$
c_{g_1,g_2}^{-1}g_1(\underline{c_{g_2,h}c_{g_2h,g_2^{-1}}})c_{g_1,g_2hg_2^{-1}}\stackrel{by (\ref{cocycle})}=c_{g_1,g_2}^{-1}g_1(g_2(c_{h,g_2^{-1}})c_{g_2,hg_2^{-1}})c_{g_1,g_2hg_2^{-1}}
$$
$$
=c_{g_1,g_2}^{-1}c_{g_1,g_2}g_1g_2(c_{h,g_2^{-1}})c_{g_1,g_2}^{-1}g_1(c_{g_2,hg_2^{-1}})c_{g_1,g_2hg_2^{-1}}=g_1g_2(c_{h,g_2^{-1}})c_{g_1,g_2}^{-1}\underline{g_1(c_{g_2,hg_2^{-1}})c_{g_1,g_2hg_2^{-1}}}
$$
$$
\stackrel{by (\ref{cocycle})}=g_1g_2(c_{h,g_2^{-1}})c_{g_1,g_2}^{-1}c_{g_1,g_2}c_{g_1g_2,hg_2^{-1}}=g_1g_2(c_{h,g_2^{-1}})c_{g_1g_2,hg_2^{-1}}.
$$
Thus, (\ref{6}) is equivalent
$$
g_1g_2(c_{h,g_2^{-1}})c_{g_1g_2,hg_2^{-1}}=c_{g_1g_2,h}c_{g_1g_2h,g_2^{-1}}
$$
which is nothing but (\ref{cocycle}) written in a different way.
The first part of the proposition is proved.

To prove the second part, we need to show that for $c\in\Ch(\cC)_g$
$$
c_e^{-1}c_{g,g^{-1}}^{-1}g(c)c_{g,g}c_{g^2,g^{-1}}=c
$$
modulo $\cC_{g}$. One has
$$
c_e^{-1}c_{g,g^{-1}}^{-1}g(c)c_{g,g}c_{g^2,g^{-1}}\equiv g(c)c_{g,g}c_{g^2,g^{-1}}g(c_e^{-1}c_{g,g^{-1}}^{-1})=g(c)c_{g,g}c_{g^2,g^{-1}}g(c_e^{-1}c_{g,g^{-1}}^{-1})
$$
$$
=g(c)c_{g,g}c_{g^2,g^{-1}}c_{g,e}^{-1}g(c_{g,g^{-1}}^{-1})\stackrel{by\,(\ref{cocycle})}{=}g(c)=1\cdot g(c)\equiv c\cdot1=c.
$$
The proposition is proved completely.

\section{Proof of Propositions \ref{alg} and \ref{coalg}}\label{a3}
\noindent {\bf Proof of Propositions \ref{alg}.}
First of all, let us point out the following property of $\xi$ which follows immediately from its definition:
for any $c\in\cC$, we have
\begin{equation}\label{inter}
c\xi'_i\otimes \xi''_i=\xi'_i\otimes \xi''_ic, \quad \xi'_ic\otimes \xi''_i =\xi'_i\otimes c\xi''_i.
\end{equation}

To prove that $m_{g,h}$ descends to a well-defined map from
$\Ch(\cC)_g\otimes\Ch(\cC)_h$ to $\Ch(\cC)_{gh}$ we need to show
that $m_{g,h}(\cC_g\otimes\cC)\subset\cC_{gh}$ and
$m_{g,h}(\cC\otimes\cC_h)\subset\cC_{gh}$.

To prove the first inclusion, let us substitute $c_1c_2-c_2g(c_1)$ for $c'$ in the expression
$\xi_i'c'g(\xi_i''c'')c_{g,h}$:
$$
\xi_i'(c_1c_2-c_2g(c_1))g(\xi_i''c'')c_{g,h}=\xi_i'c_1c_2g(\xi_i''c'')c_{g,h}-\xi_i'c_2g(c_1)g(\xi_i''c'')c_{g,h}
$$
$$
=\xi_i'c_1c_2g(\xi_i''c'')c_{g,h}-\underline{\xi_i'}c_2g(\underline{c_1\xi_i''}c'')c_{g,h}\stackrel{by (\ref{inter})}{=}\xi_i'c_1c_2g(\xi_i''c'')c_{g,h}-\xi_i'c_1c_2g(\xi_i''c'')c_{g,h}=0.
$$

Now let us substitute $c_1c_2-c_2h(c_1)$ for $c''$ in
$\xi_i'c'g(\xi_i''c'')c_{g,h}$ and show that the result belongs to $\cC_{gh}$:
$$
\xi_i'c'g(\xi_i''(c_1c_2-c_2h(c_1)))c_{g,h}=\underline{\xi_i'}c'g(\underline{\xi_i''c_1c_2})c_{g,h}-\underline{\xi_i'}c'g(\underline{\xi_i''c_2}h(c_1))c_{g,h}
$$
$$
\stackrel{by (\ref{inter})}{=}c_1c_2\xi_i'c'g(\xi_i'')c_{g,h}-c_2\xi_i'c'g(\xi_i''h(c_1))c_{g,h}=c_1c_2\xi_i'c'g(\xi_i'')c_{g,h}-c_2\xi_i'c'g(\xi_i'')g(h(c_1))c_{g,h}
$$
$$
=c_1c_2\xi_i'c'g(\xi_i'')c_{g,h}-c_2\xi_i'c'g(\xi_i'')c_{g,h}gh(c_1)c_{g,h}^{-1}c_{g,h}=c_1\left(c_2\xi_i'c'g(\xi_i'')c_{g,h}\right)-\left(c_2\xi_i'c'g(\xi_i'')c_{g,h}\right)gh(c_1).
$$
The first part of Proposition \ref{alg} is proved.

Now let us prove associativity of the maps $m_{g,h}$. We need to show that for any $g,h,k\in G$ and $c',c'',c'''\in\cC$
$$
\xi_j'\xi_i'c'g(\xi_i''c'')c_{g,h}gh(\xi_j''c''')c_{gh,k}=\xi_n'c'g(\xi_n''\xi_m'c''h(\xi_m''c''')c_{h,k})c_{g,hk}
$$
modulo $\cC_{ghk}$. In fact, we will show that the equality holds in $\cC$. Indeed, the left-hand side equals
$$
\xi_j'\xi_i'c'g(\xi_i''c'')g(h(\xi_j''c'''))c_{g,h}c_{gh,k},
$$
and the right-hand side equals
$$
\xi_n'c'g(\xi_n''\xi_m'c''h(\xi_m''c'''))g(c_{h,k})c_{g,hk}
$$
so by (\ref{cocycle}) we need to prove that
$$
\xi_j'\xi_i'c'g(\xi_i''c'')g(h(\xi_j''c'''))=\xi_n'c'g(\xi_n''\xi_m'c''h(\xi_m''c''')),
$$
or, equivalently,
$$
\xi_j'\xi_i'c'g(\xi_i''c''h(\xi_j''c'''))=\xi_n'c'g(\xi_n''\xi_m'c''h(\xi_m''c''')).
$$
The latter equality follows for (\ref{inter}), applied to $\xi_m'$ in the right-hand side. Proposition \ref{alg} is proved completely.
\medskip

\noindent{\bf Proof of Propositions \ref{coalg}.}
To prove that $\Delta_{g,h}$ descends to a map from $\Ch(\cC)_{gh}$ to $\Ch(\cC)_g\otimes\Ch(\cC)_h$, we need to show that $\Delta_{g,h}(\cC_{gh})\subset \cC_g\otimes\cC+\cC\otimes\cC_h$.
Let us substitute $c_1c_2-c_2gh(c_1)$ for $c$ in $c\,c_{g,h}^{-1}g(\xi_i')\otimes \xi_i''$:
$$
(c_1c_2-c_2gh(c_1))c_{g,h}^{-1}g(\xi_i')\otimes \xi_i''=c_1c_2c_{g,h}^{-1}g(\xi_i')\otimes \xi_i''-c_2gh(c_1)c_{g,h}^{-1}g(\xi_i')\otimes \xi_i''
$$
$$
=c_1c_2c_{g,h}^{-1}g(\xi_i')\otimes \xi_i''-c_2c_{g,h}^{-1}\underline{c_{g,h}gh(c_1)c_{g,h}^{-1}}g(\xi_i')\otimes \xi_i''
$$
$$
=c_1c_2c_{g,h}^{-1}g(\xi_i')\otimes \xi_i''-c_2c_{g,h}^{-1}g(h(c_1))g(\xi_i')\otimes \xi_i''
$$
$$
=c_1c_2c_{g,h}^{-1}g(\xi_i')\otimes \xi_i''-c_2c_{g,h}^{-1}g(h(c_1)\xi_i')\otimes \xi_i''.
$$
Modulo $\cC_{g}\otimes\cC$, the latter expression equals
$$
\xi_i'c_1c_2c_{g,h}^{-1}\otimes \xi_i''-h(c_1)\xi_i'c_2c_{g,h}^{-1}\otimes \xi_i''.
$$
Now by (\ref{inter})
$$
\xi_i'c_1c_2c_{g,h}^{-1}\otimes \xi_i''-h(c_1)\xi_i'c_2c_{g,h}^{-1}\otimes \xi_i''=
\xi_i'c_2c_{g,h}^{-1}\otimes c_1\xi_i''-\xi_i'c_2c_{g,h}^{-1}\otimes \xi_i''h(c_1)\in\cC\otimes\cC_h.
$$
The first part of the proposition is proved.

To prove coassociativity, observe that
$$
(\Delta_{g,h}\otimes1)\Delta_{gh,k}(c)=c\,c_{gh,k}^{-1}gh(\xi_i')c_{g,h}^{-1}g(\xi_j')\otimes \xi_j''\otimes \xi_i'',
$$
$$
(1\otimes\Delta_{h,k})\Delta_{g,hk}(c)=c\,c_{g,hk}^{-1}g(\xi_j')\otimes \xi_j''c_{h,k}^{-1}h(\xi_i')\otimes \xi_i''\stackrel{by (\ref{inter})}{=}c\,c_{g,hk}^{-1}g(c_{h,k}^{-1}h(\xi_i')\xi_j')\otimes \xi_j''\otimes \xi_i''.
$$
Therefore, it is enough to show that
$$
c_{gh,k}^{-1}gh(\xi_i')c_{g,h}^{-1}g(\xi_j')=c_{g,hk}^{-1}g(c_{h,k}^{-1}h(\xi_i')\xi_j')
$$
or, equivalently,
$$
c_{gh,k}^{-1}gh(\xi_i')c_{g,h}^{-1}g(\xi_j')=c_{g,hk}^{-1}g(c_{h,k}^{-1})c_{g,h}gh(\xi_i')c_{g,h}^{-1}g(\xi_j').
$$
This is an immediate consequence of (\ref{cocycle}).
\medskip
\section{Proof of Theorem \ref{main}}\label{a4}
We need to prove the $G$-invariance of $\theta_e$ and properties
(1), (3), (4), (5), (7), (8), and (9) from Definitions \ref{def1}
and \ref{def2} (properties (2) and (6) are satisfied by construction
of the multiplication and the comultiplication). As before, to make
following our computations easier, we will sometimes underline those
places in formulas that are about to be changed.

\medskip
\noindent{\bf Proof of the $G$-invariance of $\theta_e$.} For any
$c\in\Ch(\cC)_e$ and $g\in G$
$$
\theta_e(T_e(g)(c))=\theta(c_e^{-1}c_{g,g^{-1}}^{-1}g(c)c_{g,e}c_{g,g^{-1}}c_e)=\theta(c_{g,g^{-1}}^{-1}g(c)c_{g,e}c_{g,g^{-1}})
$$
$$
=\theta(g(c)c_{g,e})=\theta(g(c)g(c_{e}))=\theta(g(cc_{e}))=\theta(cc_{e})=\theta_e(c).
$$

\medskip
\noindent{\bf Proof of property (1).}
Let $c'\in\Ch(\cC)_h$, $c''\in\Ch(\cC)_k$, and $g\in G$. We need to show that
$$
m_{ghg^{-1},gkg^{-1}}(\cT_h(g)(c')\otimes\cT_k(g)(c''))=\cT_{hk}(g)(m_{h,k}(c'\otimes c'')).
$$
We have
$$
m_{ghg^{-1},gkg^{-1}}(\cT_h(g)(c')\otimes\cT_k(g)(c''))
$$
$$
=m_{ghg^{-1},gkg^{-1}}((c_e^{-1}c_{g,g^{-1}}^{-1}g(c')c_{g,h}c_{gh,g^{-1}})\otimes (c_e^{-1}c_{g,g^{-1}}^{-1}g(c'')c_{g,k}c_{gk,g^{-1}}))
$$
$$
=\xi_i'c_e^{-1}c_{g,g^{-1}}^{-1}g(c')c_{g,h}c_{gh,g^{-1}}ghg^{-1}(\xi_i''c_e^{-1}c_{g,g^{-1}}^{-1}g(c'')c_{g,k}c_{gk,g^{-1}})c_{ghg^{-1},gkg^{-1}}
$$
$$
=\xi_i'c_e^{-1}c_{g,g^{-1}}^{-1}g(c')c_{g,h}c_{gh,g^{-1}}ghg^{-1}(\xi_i'')ghg^{-1}(c_e^{-1}c_{g,g^{-1}}^{-1})c_{ghg^{-1},g}gh(c'')\times
$$
$$
\times \underline{c_{ghg^{-1},g}^{-1}ghg^{-1}(c_{g,k})ghg^{-1}(c_{gk,g^{-1}})c_{ghg^{-1},gkg^{-1}}}
$$
$$
\stackrel{by (\ref{cocycle})}{=}\xi_i'c_e^{-1}c_{g,g^{-1}}^{-1}g(c')c_{g,h}c_{gh,g^{-1}}ghg^{-1}(\xi_i'')\underline{ghg^{-1}(c_e^{-1}c_{g,g^{-1}}^{-1})c_{ghg^{-1},g}}gh(c'')c_{gh,k}c_{ghk,g^{-1}}
$$
$$
\stackrel{by (\ref{cocycle})}{=}\xi_i'c_e^{-1}c_{g,g^{-1}}^{-1}g(c')c_{g,h}\underline{c_{gh,g^{-1}}ghg^{-1}(\xi_i'')c_{gh,g^{-1}}^{-1}}gh(c'')c_{gh,k}c_{ghk,g^{-1}}
$$
$$
=\xi_i'c_e^{-1}c_{g,g^{-1}}^{-1}g(c')c_{g,h}gh(g^{-1}(\xi_i''))gh(c'')c_{gh,k}c_{ghk,g^{-1}}
$$
$$
=\xi_i'c_e^{-1}c_{g,g^{-1}}^{-1}g(c')c_{g,h}gh(g^{-1}(\xi_i'')c'')c_{gh,k}c_{ghk,g^{-1}}.
$$
On the other hand,
$$
\cT_{hk}(g)(m_{h,k}(c'\otimes c''))=\cT_{hk}(g)(\xi_i'c'h(\xi_i''c'')c_{h,k})
=c_e^{-1}c_{g,g^{-1}}^{-1}g(\xi_i'c'h(\xi_i''c'')c_{h,k})c_{g,hk}c_{ghk,g^{-1}}
$$
$$
=c_e^{-1}c_{g,g^{-1}}^{-1}g(\xi_i')g(c')c_{g,h}gh(\xi_i''c'')\underline{c_{g,h}^{-1}g(c_{h,k})c_{g,hk}}c_{ghk,g^{-1}}
$$
$$
\stackrel{by (\ref{cocycle})}{=}c_e^{-1}c_{g,g^{-1}}^{-1}g(\xi_i')g(c')c_{g,h}gh(\xi_i''c'')c_{gh,k}c_{ghk,g^{-1}}.
$$
So we have to prove that
$$
\xi_i'c_e^{-1}c_{g,g^{-1}}^{-1}g(c')c_{g,h}gh(g^{-1}(\xi_i'')c'')=c_e^{-1}c_{g,g^{-1}}^{-1}g(\xi_i')g(c')c_{g,h}gh(\xi_i''c'').
$$
The latter equality is an immediate consequence of the following
fact which is just another way to express the $G$-invariance of the
trace $\theta$: for any $g\in G$,
\begin{equation}\label{Ginv}
\xi_i'\otimes g^{-1}(\xi_i'')=Ad(c_e^{-1})Ad(c_{g,g^{-1}}^{-1})g(\xi_i')\otimes \xi_i''.
\end{equation}

\medskip
\noindent{\bf Proof of property (3).}
Let $c'\in\Ch(\cC)_g$, $c''\in\Ch(\cC)_h$, and $g\in G$. We need to show that
$$
m_{g,h}(c'\otimes c'')=m_{ghg^{-1},g}(\cT_h(g)(c'')\otimes c')
$$
modulo $\cC_{gh}$. We have
$$
m_{ghg^{-1},g}(\cT_h(g)(c'')\otimes c')=m_{ghg^{-1},g}(c_e^{-1}c_{g,g^{-1}}^{-1}g(c'')c_{g,h}c_{gh,g^{-1}}\otimes c')
$$
$$
=\xi_i'c_e^{-1}c_{g,g^{-1}}^{-1}g(c'')c_{g,h}c_{gh,g^{-1}}ghg^{-1}(\xi_i''c')c_{ghg^{-1},g}
$$
$$
=\xi_i'c_e^{-1}c_{g,g^{-1}}^{-1}g(c'')c_{g,h}\underline{c_{gh,g^{-1}}ghg^{-1}(\xi_i''c')c_{gh,g^{-1}}^{-1}}c_{gh,g^{-1}}c_{ghg^{-1},g}
$$
$$
=\xi_i'c_e^{-1}c_{g,g^{-1}}^{-1}g(c'')c_{g,h}gh(g^{-1}(\xi_i''c'))\underline{c_{gh,g^{-1}}c_{ghg^{-1},g}}
$$
$$
=\xi_i'c_e^{-1}c_{g,g^{-1}}^{-1}g(c'')c_{g,h}gh(g^{-1}(\xi_i''c'))gh(c_{g^{-1},g})\underline{c_{gh,e}}
$$
$$
=\xi_i'c_e^{-1}c_{g,g^{-1}}^{-1}g(c'')c_{g,h}\underline{gh(g^{-1}(\xi_i''c'))gh(c_{g^{-1},g})gh(c_{e})}
$$
$$
=\xi_i'c_e^{-1}c_{g,g^{-1}}^{-1}g(c'')c_{g,h}gh(g^{-1}(\xi_i''c')c_{g^{-1},g}c_{e}).
$$
Modulo $\cC_{gh}$, the latter expression equals
$$
g^{-1}(\xi_i''c')c_{g^{-1},g}c_{e}\xi_i'c_e^{-1}c_{g,g^{-1}}^{-1}g(c'')c_{g,h}.
$$
Then
$$
g^{-1}(\xi_i''c')c_{g^{-1},g}c_{e}\xi_i'c_e^{-1}c_{g,g^{-1}}^{-1}g(c'')c_{g,h}\stackrel{by (\ref{inter})}{=}
g^{-1}(\xi_i'')c_{g^{-1},g}c_{e}c'\xi_i'c_e^{-1}c_{g,g^{-1}}^{-1}g(c'')c_{g,h}
$$
$$
\stackrel{by (\ref{Ginv})}{=}\xi_i''c_{g^{-1},g}c_{e}c'c_e^{-1}c_{g,g^{-1}}^{-1}g(\xi_i')c_{g,g^{-1}}c_ec_e^{-1}c_{g,g^{-1}}^{-1}g(c'')c_{g,h}=\xi_i''c_{g^{-1},g}c_{e}c'c_e^{-1}c_{g,g^{-1}}^{-1}g(\xi_i')g(c'')c_{g,h}
$$
$$
=\xi_i''c_{g^{-1},g}c_{e}c'c_e^{-1}c_{g,g^{-1}}^{-1}g(\xi_i'c'')c_{g,h}\stackrel{by (\ref{inter})}{=}\xi_i''c'c_e^{-1}c_{g,g^{-1}}^{-1}g(c_{g^{-1},g}c_{e}\xi_i'c'')c_{g,h}
$$
$$
{=}\xi_i''c'c_e^{-1}c_{g,g^{-1}}^{-1}g(c_{g^{-1},g}c_{e})g(\xi_i'c'')c_{g,h}\stackrel{by (\ref{cocycle})}{=}\xi_i''c'g(\xi_i'c'')c_{g,h}.
$$
The latter expression equals $m_{g,h}(c'\otimes c'')$ since $\xi$ is symmetric.

\medskip
\noindent{\bf Proof of property (5).}
We need to show that for any $g,h,k\in G$ and $c\in\Ch(\cC)_{hk}$
$$
\Delta_{ghg^{-1},gkg^{-1}}(\cT_{hk}(g)(c))=(\cT_{h}(g)\otimes\cT_{k}(g))\Delta_{h,k}(c).
$$
We have
$$
(\cT_{h}(g)\otimes\cT_{k}(g))\Delta_{h,k}(c)=(\cT_{h}(g)\otimes\cT_{k}(g))(c\,c_{h,k}^{-1}h(\xi_i')\otimes \xi_i'')
$$
$$
=c_e^{-1}c_{g,g^{-1}}^{-1}g(c\,c_{h,k}^{-1}h(\xi_i'))c_{g,h}c_{gh,g^{-1}}\otimes c_e^{-1}c_{g,g^{-1}}^{-1}\underline{g(\xi_i'')c_{g,k}}c_{gk,g^{-1}}
$$
$$
=c_e^{-1}c_{g,g^{-1}}^{-1}g(c\,c_{h,k}^{-1}h(\underline{\xi_i'}))c_{g,h}c_{gh,g^{-1}}\otimes \underline{c_e^{-1}c_{g,g^{-1}}^{-1}g(\xi_i'')c_{g,g^{-1}}c_e}c_e^{-1}c_{g,g^{-1}}^{-1}c_{g,k}c_{gk,g^{-1}}
$$
$$
\stackrel{by (\ref{Ginv})}{=}c_e^{-1}c_{g,g^{-1}}^{-1}\underline{g(c\,c_{h,k}^{-1}h(g^{-1}(\xi_i')))}c_{g,h}c_{gh,g^{-1}}\otimes \xi_i''c_e^{-1}c_{g,g^{-1}}^{-1}c_{g,k}c_{gk,g^{-1}}
$$
$$
=c_e^{-1}c_{g,g^{-1}}^{-1}g(c)g(c_{h,k}^{-1})\underline{g(h(g^{-1}(\xi_i')))}c_{g,h}c_{gh,g^{-1}}\otimes \xi_i''c_e^{-1}c_{g,g^{-1}}^{-1}c_{g,k}c_{gk,g^{-1}}
$$
$$
=c_e^{-1}c_{g,g^{-1}}^{-1}g(c)g(c_{h,k}^{-1})\underline{g(c_{h,g^{-1}}hg^{-1}(\xi_i')c_{h,g^{-1}}^{-1})}c_{g,h}c_{gh,g^{-1}}\otimes \xi_i''c_e^{-1}c_{g,g^{-1}}^{-1}c_{g,k}c_{gk,g^{-1}}
$$
$$
=c_e^{-1}c_{g,g^{-1}}^{-1}g(c)g(c_{h,k}^{-1})g(c_{h,g^{-1}})c_{g,hg^{-1}}ghg^{-1}(\xi_i')\underline{c_{g,hg^{-1}}^{-1}g(c_{h,g^{-1}}^{-1})c_{g,h}c_{gh,g^{-1}}}\otimes \xi_i''c_e^{-1}c_{g,g^{-1}}^{-1}c_{g,k}c_{gk,g^{-1}}
$$
$$
\stackrel{by (\ref{cocycle})}{=}c_e^{-1}c_{g,g^{-1}}^{-1}g(c)g(c_{h,k}^{-1})g(c_{h,g^{-1}})c_{g,hg^{-1}}ghg^{-1}(\underline{\xi_i'})\otimes \underline{\xi_i''c_e^{-1}c_{g,g^{-1}}^{-1}c_{g,k}c_{gk,g^{-1}}}
$$
$$
\stackrel{by (\ref{inter})}{=}c_e^{-1}c_{g,g^{-1}}^{-1}g(c)g(c_{h,k}^{-1})\underline{g(c_{h,g^{-1}})c_{g,hg^{-1}}}
\underline{ghg^{-1}(c_e^{-1}c_{g,g^{-1}}^{-1}}c_{g,k}c_{gk,g^{-1}})
ghg^{-1}(\xi_i')\otimes \xi_i''
$$
$$
\stackrel{by (\ref{cocycle})}{=}c_e^{-1}c_{g,g^{-1}}^{-1}g(c)g(c_{h,k}^{-1})c_{g,h}c_{gh,g^{-1}}
c_{gh,g^{-1}}^{-1}c_{ghg^{-1},g}^{-1}ghg^{-1}(c_{g,k}c_{gk,g^{-1}})
ghg^{-1}(\xi_i')\otimes \xi_i''
$$
$$
{=}c_e^{-1}c_{g,g^{-1}}^{-1}g(c)g(c_{h,k}^{-1})c_{g,h}c_{ghg^{-1},g}^{-1}ghg^{-1}(c_{g,k}c_{gk,g^{-1}})
ghg^{-1}(\xi_i')\otimes \xi_i''.
$$
On the other hand,
$$
\Delta_{ghg^{-1},gkg^{-1}}(\cT_{hk}(g)(c))=\Delta_{ghg^{-1},gkg^{-1}}(c_e^{-1}c_{g,g^{-1}}^{-1}g(c)c_{g,hk}c_{ghk,g^{-1}})
$$
$$
=c_e^{-1}c_{g,g^{-1}}^{-1}g(c)c_{g,hk}c_{ghk,g^{-1}}c_{ghg^{-1},gkg^{-1}}^{-1}ghg^{-1}(\xi_i')\otimes \xi_i''.
$$
Thus, we need to show that
$$
g(c_{h,k}^{-1})c_{g,h}c_{ghg^{-1},g}^{-1}ghg^{-1}(c_{g,k}c_{gk,g^{-1}})=c_{g,hk}c_{ghk,g^{-1}}c_{ghg^{-1},gkg^{-1}}^{-1}
$$
or, equivalently,
$$
c_{g,h}c_{ghg^{-1},g}^{-1}ghg^{-1}(c_{g,k})ghg^{-1}(c_{gk,g^{-1}})c_{ghg^{-1},gkg^{-1}}=g(c_{h,k})c_{g,hk}c_{ghk,g^{-1}}.
$$
Let us transform both hand sides of the latter equality
simultaneously. First of all, (\ref{cocycle}), applied to the
following parts of the expressions
$$
c_{g,h}c_{ghg^{-1},g}^{-1}ghg^{-1}(c_{g,k})\underline{ghg^{-1}(c_{gk,g^{-1}})c_{ghg^{-1},gkg^{-1}}}=\underline{g(c_{h,k})c_{g,hk}}c_{ghk,g^{-1}},
$$
gives us
$$
c_{g,h}c_{ghg^{-1},g}^{-1}ghg^{-1}(c_{g,k})c_{ghg^{-1},gk}c_{ghk,g^{-1}}=c_{g,h}c_{gh,k}c_{ghk,g^{-1}},
$$
which is equivalent to
$$
c_{ghg^{-1},g}^{-1}ghg^{-1}(c_{g,k})c_{ghg^{-1},gk}=c_{gh,k}.
$$
The latter is equivalent to (\ref{cocycle}).

\medskip
\noindent{\bf Proof of property (7).}
Let $c\in\Ch(\cC)_{gh}$. We need to show that
$$
\Delta_{g,h}(c)=\sigma (1\otimes \cT_{h^{-1}gh}(h))(\Delta_{h,h^{-1}gh}(c))
$$
in $\Ch(\cC)_{g}\otimes\Ch(\cC)_{h}$.
We have
$$
(1\otimes \cT_{h^{-1}gh}(h))(\Delta_{h,h^{-1}gh}(c))=(1\otimes \cT_{h^{-1}gh}(h))(c\,c_{h,h^{-1}gh}^{-1}h(\xi_i')\otimes \xi_i'')
$$
$$
=(1\otimes \cT_{h^{-1}gh}(h))(c\,c_{h,h^{-1}gh}^{-1}h(\xi_i')\otimes \xi_i'')
=c\,c_{h,h^{-1}gh}^{-1}h(\xi_i')\otimes c_e^{-1}c_{h,h^{-1}}^{-1}h(\xi_i'')c_{h,h^{-1}gh}c_{gh,h^{-1}}.
$$
Now observe that (\ref{Ginv}) is equivalent
$$
h(\xi_i')\otimes h(\xi_i'')=Ad(c_{h,h^{-1}})Ad(c_e)(\xi_i')\otimes Ad(c_{h,h^{-1}})Ad(c_e)(\xi_i'')\stackrel{by(\ref{inter})}=\xi_i'\otimes \xi_i''.
$$
Therefore,
$$
(1\otimes \cT_{h^{-1}gh}(h))(\Delta_{h,h^{-1}gh}(c))
=c\,c_{h,h^{-1}gh}^{-1}\xi_i'\otimes c_e^{-1}c_{h,h^{-1}}^{-1}\xi_i''c_{h,h^{-1}gh}c_{gh,h^{-1}}
$$
$$
\stackrel{by (\ref{inter})}=\xi_i'\otimes c_e^{-1}c_{h,h^{-1}}^{-1}\xi_i''c\,c_{h,h^{-1}gh}^{-1}c_{h,h^{-1}gh}c_{gh,h^{-1}}.
$$
Thus,
$$
\sigma(1\otimes \cT_{h^{-1}gh}(h))(\Delta_{h,h^{-1}gh}(c))=c_e^{-1}c_{h,h^{-1}}^{-1}\xi_i'c\,c_{h,h^{-1}gh}^{-1}c_{h,h^{-1}gh}c_{gh,h^{-1}}\otimes \xi_i''.
$$
Modulo $\cC_g\otimes\cC$, the latter expression equals
$$
c\,c_{h,h^{-1}gh}^{-1}c_{h,h^{-1}gh}c_{gh,h^{-1}}g(c_e^{-1}c_{h,h^{-1}}^{-1}\xi_i')\otimes \xi_i''=c\,c_{gh,h^{-1}}g(c_e^{-1}c_{h,h^{-1}}^{-1}\xi_i')\otimes \xi_i''.
$$
To complete the proof, observe that
$$
c\,c_{gh,h^{-1}}\underline{g(c_e^{-1}c_{h,h^{-1}}^{-1}\xi_i')}\otimes \xi_i''
\stackrel{by (\ref{cocycle})}=c\,c_{gh,h^{-1}}c_{gh,h^{-1}}^{-1}c_{g,h}^{-1}g(\xi_i')\otimes \xi_i''
=c\,c_{g,h}^{-1}g(\xi_i')\otimes \xi_i''=\Delta_{g,h}(c).
$$

\medskip
\noindent{\bf Proof of property (8).}
Let $g\in G$ and $c\in\Ch(\cC)_g$. Then
$$
(\theta_e\otimes1)\Delta_{e,g}(c)=(\theta_e\otimes1)(c\,c_{e,g}^{-1}e(\xi_i')\otimes \xi_i'')\stackrel{by (\ref{cocycle}) and (\ref{theta})}=\theta(c\,c_{e}^{-1}e(\xi_i')c_e)\xi_i''=\theta(c\,\xi_i')\xi_i''=c.
$$
Similarly,
$$
(1\otimes\theta_e)\Delta_{g,e}(c)=(1\otimes\theta_e)(c\,c_{g,e}^{-1}g(\xi_i')\otimes \xi_i'')
=c\,c_{g,e}^{-1}g(\xi_i')\theta_e(\xi_i'')=c\,c_{g,e}^{-1}g(\xi_i')\theta(\xi_i''c_e)
$$
$$
=c\,c_{g,e}^{-1}g(c_e)=c.
$$

\medskip
\noindent{\bf Proof of property (9).}
Let $g,h,k\in G$, $c'\in\Ch(\cC)_g$, $c''\in\Ch(\cC)_{hk}$. Firstly, we need to prove that
$\Delta$ is a morphism of left $\Ch(\cC)$-modules, i.e.
$$
\Delta_{gh,k}(m_{g,hk}(c'\otimes c''))=(m_{g,h}\otimes1)(c'\otimes \Delta_{h,k}(c''))
$$
We have
$$
\Delta_{gh,k}(m_{g,hk}(c'\otimes c''))=\Delta_{gh,k}(\xi_i'c'g(\xi_i''c'')c_{g,hk})
=\xi_i'c'g(\xi_i''c'')c_{g,hk}c_{gh,k}^{-1}gh(\xi_j')\otimes \xi_j''.
$$
On the other hand,
$$
(m_{g,h}\otimes1)(c'\otimes \Delta_{h,k}(c''))=(m_{g,h}\otimes1)(c'\otimes c''\,c_{h,k}^{-1}h(\xi_j')\otimes \xi_j'')=
\xi_i'c'g(\xi_i''c''c_{h,k}^{-1}h(\xi_j'))c_{g,h}\otimes \xi_j''
$$
$$
=\xi_i'c'g(\xi_i''c'')g(c_{h,k}^{-1})c_{g,h}gh(\xi_j')c_{g,h}^{-1}c_{g,h}\otimes \xi_j''=\xi_i'c'g(\xi_i''c'')\underline{g(c_{h,k}^{-1})c_{g,h}}gh(\xi_j')\otimes \xi_j''
$$
$$
\stackrel{by (\ref{cocycle})}=\xi_i'c'g(\xi_i''c'')c_{g,hk}c_{gh,k}^{-1}gh(\xi_j')\otimes \xi_j''.
$$

Secondly, we need to prove that $\Delta$ is a morphism of right $\Ch(\cC)$-modules, i.e.
\begin{equation}\label{rm}
\Delta_{h,kg}(m_{hk,g}(c''\otimes c'))=(1\otimes m_{k,g})(\Delta_{h,k}(c'')\otimes c')
\end{equation}
in $\Ch(\cC)_h\otimes \Ch(\cC)_{kg}$. We have
$$
\Delta_{h,kg}(m_{hk,g}(c''\otimes c'))=\Delta_{h,kg}(\xi_j'c''hk(\xi_j''c')c_{hk,g})=\xi_j'c''hk(\xi_j''c')c_{hk,g}c_{h,kg}^{-1}h(\xi_i')\otimes \xi_i''.
$$
On the other hand,
\begin{equation}\label{7}
(1\otimes m_{k,g})(\Delta_{h,k}(c'')\otimes c')=(1\otimes m_{k,g})(c''c_{h,k}^{-1}h(\xi_i')\otimes \xi_i''\otimes c')=c''c_{h,k}^{-1}h(\underline{\xi_i'})\otimes\underline{\xi_j'\xi_i''k(\xi_j''c')c_{k,g}}
\end{equation}
$$
\stackrel{by (\ref{inter})}=c''c_{h,k}^{-1}h(k(\xi_j''c')c_{k,g}\xi_i'\xi_j')\otimes\xi_i''=c''hk(\xi_j''c')c_{h,k}^{-1}h(c_{k,g}\xi_i'\xi_j')\otimes\xi_i''
$$
$$
=c''hk(\xi_j''c')\underline{c_{h,k}^{-1}h(c_{k,g})}h(\xi_i')h(\xi_j')\otimes\xi_i''\stackrel{by (\ref{cocycle})}=c''hk(\xi_j''c')c_{hk,g}c_{h,kg}^{-1}h(\xi_i')h(\xi_j')\otimes\xi_i''.
$$
Thus, we see that
$$
\Delta_{h,kg}(m_{hk,g}(c''\otimes c'))-(1\otimes m_{k,g})(\Delta_{h,k}(c'')\otimes c')
$$
$$
=\xi_j'\left(c''hk(\xi_j''c')c_{hk,g}c_{h,kg}^{-1}h(\xi_i')\right)\otimes \xi_i''-\left(c''hk(\xi_j''c')c_{hk,g}c_{h,kg}^{-1}h(\xi_i')\right)h(\xi_j')\otimes\xi_i''\in\cC_h\otimes\cC.
$$

\medskip
\noindent{\bf Proof of property (4) (the torus axiom).}
First of all, let us reformulate the torus axiom. Let $C$ be an object satisfying all the axioms of a weak crossed $G$-algebra except for the torus axiom. For $g,h\in G$, let $m_{g,h}$ stand for the multiplication map $C_g\otimes C_h\to C_{gh}$.

Fix $g_1,g_2\in G$, $c\in C_{g_1}$, $c'\in C_{g_2}$. Then
$$
(1\otimes\theta_e)(1\otimes m_{g_2^{-1},g_2})(\Delta_{g_1g_2,g_2^{-1}}(c)\otimes c')\stackrel{by \,(8)}=(1\otimes\theta_e)(\Delta_{g_1g_2,e}(m_{g_1,g_2}(c\otimes c')))
\stackrel{by\,(7)}=m_{g_1,g_2}(c\otimes c').
$$
If we set
$
\Delta_{g_1g_2,g_2^{-1}}(c)=\sum c^{(1)}\otimes c^{(2)}
$
then the above computation means that under the canonical isomorphism
$$
\Hom(C_{g_2},C_{g_1g_2})\cong C_{g_1g_2}\otimes C_{g_2}^*,
$$
the operator $m_{g_1,g_2}(c\otimes \bullet)$ corresponds to the element
$$
c^{(1)}\otimes \theta_e(m_{g_2^{-1},g_2}(c^{(2)}\otimes \bullet)).
$$
Thus, if $T:C_{g_1g_2}\to C_{g_2}$ is a linear operator then
$$
{\Tr}_{C_{g_1g_2}}(m_{g_1,g_2}(c\otimes \bullet)\cdot T)={\Tr}_{C_{g_1g_2}}(c^{(1)}\otimes \theta_e(m_{g_2^{-1},g_2}(c^{(2)}\otimes T(\bullet)))=\theta_e(m_{g_2^{-1},g_2}(c^{(2)}\otimes T(c^{(1)})))
$$
\begin{equation}\label{11}
=\theta_e(m_{g_2^{-1},g_2}\cdot\sigma\cdot(T\otimes 1)\cdot\Delta_{g_1g_2,g_2^{-1}}(c))
\end{equation}
and, similarly,
$$
{\Tr}_{C_{g_2}}(T\cdot m_{g_1,g_2}(c\otimes \bullet))={\Tr}_{C_{g_2}}(T(c^{(1)})\otimes \theta_e(m_{g_2^{-1},g_2}(c^{(2)}\otimes \bullet)))=\theta_e(m_{g_2^{-1},g_2}(c^{(2)}\otimes T(c^{(1)})))
$$
\begin{equation}\label{22}
=\theta_e(m_{g_2^{-1},g_2}\cdot\sigma\cdot(T\otimes 1)\cdot\Delta_{g_1g_2,g_2^{-1}}(c)).
\end{equation}

Now we are ready to reformulate the torus axiom. Let us fix $g,h\in G$ and $c\in C_{hgh^{-1}g^{-1}}$.
Then, by (\ref{11})
$$
{\Tr}_{C_h}(m_{hgh^{-1}g^{-1},ghg^{-1}}(c\otimes \bullet)\cdot g)=\theta_e(m_{gh^{-1}g^{-1},ghg^{-1}}\cdot\sigma\cdot(g\otimes 1)\cdot\Delta_{h,gh^{-1}g^{-1}}(c))
$$
$$
\stackrel{by\,G-inv.}=\theta_e(g^{-1}\cdot m_{gh^{-1}g^{-1},ghg^{-1}}\cdot\sigma\cdot(g\otimes 1)\cdot\Delta_{h,gh^{-1}g^{-1}}(c))
$$
$$
\stackrel{by\,(1)}=\theta_e(m_{h^{-1},h}\cdot(g^{-1}\otimes g^{-1})\cdot\sigma\cdot(g\otimes 1)\cdot\Delta_{h,gh^{-1}g^{-1}}(c))
$$
$$
=\theta_e(m_{h^{-1},h}\cdot\sigma\cdot(1\otimes g^{-1})\cdot\Delta_{h,gh^{-1}g^{-1}}(c))
$$
$$
\stackrel{by\,(3)}=\theta_e(m_{h,h^{-1}}\cdot(1\otimes g^{-1})\cdot\Delta_{h,gh^{-1}g^{-1}}(c)).
$$
On the other hand, by (\ref{22})
$$
{\Tr}_{C_{g}}(h^{-1}\cdot m_{hgh^{-1}g^{-1},g}(c\otimes\bullet))=\theta_e(m_{g^{-1},g}\cdot\sigma\cdot(h^{-1}\otimes 1)\cdot\Delta_{hgh^{-1},g^{-1}}(c))
$$
$$
\stackrel{by\,(3)}=\theta_e(m_{g,g^{-1}}\cdot(h^{-1}\otimes 1)\cdot\Delta_{hgh^{-1},g^{-1}}(c))
$$
$$
\stackrel{by\,G-inv.}=\theta_e(h\cdot m_{g,g^{-1}}\cdot(h^{-1}\otimes 1)\cdot\Delta_{hgh^{-1},g^{-1}}(c))
$$
$$
\stackrel{by\,(1)}=\theta_e(m_{hgh^{-1},hg^{-1}h^{-1}}\cdot(1\otimes h)\cdot\Delta_{hgh^{-1},g^{-1}}(c)).
$$
Therefore, the torus axiom is equivalent to the following statement: for any $g,h\in G$ and $c\in C_{hgh^{-1}g^{-1}}$
\begin{equation}\label{newta}
\theta_e(m_{h,h^{-1}}\cdot(1\otimes g^{-1})\cdot\Delta_{h,gh^{-1}g^{-1}}(c))=\theta_e(m_{hgh^{-1},hg^{-1}h^{-1}}\cdot(1\otimes h)\cdot\Delta_{hgh^{-1},g^{-1}}(c)).
\end{equation}
We will show that (\ref{newta}) is satisfied in $\Ch(\cC)$.

Consider an element $c\in\Ch(\cC)_{hgh^{-1}g^{-1}}$. Then
$$
\theta_e(m_{h,h^{-1}}\cdot(1\otimes \cT_{gh^{-1}g^{-1}}(g^{-1}))\cdot\Delta_{h,gh^{-1}g^{-1}}(c))
$$
$$
=\theta_e(m_{h,h^{-1}}\cdot(1\otimes \cT_{gh^{-1}g^{-1}}(g^{-1}))(c\,c_{h,gh^{-1}g^{-1}}^{-1}h(\xi_i')\otimes \xi_i''))
$$
$$
=\theta_e(m_{h,h^{-1}}(c\,c_{h,gh^{-1}g^{-1}}^{-1}h(\xi_i')\otimes c_e^{-1}c_{g^{-1},g}^{-1}g^{-1}(\xi_i'')c_{g^{-1},gh^{-1}g^{-1}}c_{h^{-1}g^{-1},g}))
$$
$$
=\theta_e(\xi_j'c\,c_{h,gh^{-1}g^{-1}}^{-1}h(\xi_i')h(\xi_j''c_e^{-1}c_{g^{-1},g}^{-1}g^{-1}(\xi_i'')c_{g^{-1},gh^{-1}g^{-1}}c_{h^{-1}g^{-1},g})c_{h,h^{-1}})
$$
$$
=\theta_e(\xi_j'c\,c_{h,gh^{-1}g^{-1}}^{-1}h(\xi_i')h(\xi_j'')\underline{h(c_e^{-1})h(c_{g^{-1},g}^{-1})}h(g^{-1}(\xi_i''))\underline{h(c_{g^{-1},gh^{-1}g^{-1}})h(c_{h^{-1}g^{-1},g})c_{h,h^{-1}}})
$$
$$
\stackrel{by (\ref{cocycle})}=\theta_e(\xi_j'c\,c_{h,gh^{-1}g^{-1}}^{-1}h(\xi_i')h(\xi_j'')c_{hg^{-1},g}^{-1}\underline{c_{h,g^{-1}}^{-1}h(g^{-1}(\xi_i''))c_{h,g^{-1}}}c_{hg^{-1},gh^{-1}g^{-1}}c_{g^{-1},g})
$$
$$
=\theta_e(\xi_j'c\,c_{h,gh^{-1}g^{-1}}^{-1}h(\xi_i')h(\xi_j'')c_{hg^{-1},g}^{-1}hg^{-1}(\xi_i'')c_{hg^{-1},gh^{-1}g^{-1}}c_{g^{-1},g})
$$
\begin{equation}\label{33}
=\theta(c\,c_{h,gh^{-1}g^{-1}}^{-1}h(\xi_i')h(\xi_j'')c_{hg^{-1},g}^{-1}hg^{-1}(\xi_i'')c_{hg^{-1},gh^{-1}g^{-1}}c_{g^{-1},g}c_e\xi_j').
\end{equation}
On the other hand,
$$
\theta_e(m_{hgh^{-1},hg^{-1}h^{-1}}\cdot(1\otimes \cT_{g^{-1}}(h))\cdot\Delta_{hgh^{-1},g^{-1}}(c))
$$
$$
=\theta_e(m_{hgh^{-1},hg^{-1}h^{-1}}\cdot(1\otimes \cT_{g^{-1}}(h))(c\,c_{hgh^{-1},g^{-1}}^{-1}hgh^{-1}(\xi_i')\otimes \xi_i''))
$$
$$
=\theta_e(m_{hgh^{-1},hg^{-1}h^{-1}}(c\,c_{hgh^{-1},g^{-1}}^{-1}hgh^{-1}(\xi_i')\otimes c_e^{-1}c_{h,h^{-1}}^{-1}h(\xi_i'')c_{h,g^{-1}}c_{hg^{-1},h^{-1}}))
$$
$$
=\theta_e(\xi_j'c\,c_{hgh^{-1},g^{-1}}^{-1}hgh^{-1}(\xi_i')hgh^{-1}(\xi_j''c_e^{-1}c_{h,h^{-1}}^{-1}h(\xi_i'')c_{h,g^{-1}}c_{hg^{-1},h^{-1}})c_{hgh^{-1},hg^{-1}h^{-1}})
$$
$$
=\theta_e(\xi_j'c\,c_{hgh^{-1},g^{-1}}^{-1}hgh^{-1}(\xi_i')hgh^{-1}(\xi_j'')\underline{hgh^{-1}(c_e^{-1}c_{h,h^{-1}}^{-1})}hgh^{-1}(h(\xi_i''))\times
$$
$$
\times \underline{hgh^{-1}(c_{h,g^{-1}})hgh^{-1}(c_{hg^{-1},h^{-1}})c_{hgh^{-1},hg^{-1}h^{-1}}})
$$
$$
=\theta_e(\xi_j'c\,c_{hgh^{-1},g^{-1}}^{-1}hgh^{-1}(\xi_i')hgh^{-1}(\xi_j'')c_{hg,h^{-1}}^{-1}\underline{c_{hgh^{-1},h}^{-1}hgh^{-1}(h(\xi_i''))c_{hgh^{-1},h}}c_{hg,g^{-1}}c_{h,h^{-1}})
$$
$$
=\theta_e(\xi_j'c\,c_{hgh^{-1},g^{-1}}^{-1}hgh^{-1}(\xi_i')hgh^{-1}(\xi_j'')c_{hg,h^{-1}}^{-1}hg(\xi_i'')c_{hg,g^{-1}}c_{h,h^{-1}})
$$
\begin{equation}\label{44}
=\theta(c\,c_{hgh^{-1},g^{-1}}^{-1}hgh^{-1}(\xi_i')hgh^{-1}(\xi_j'')c_{hg,h^{-1}}^{-1}hg(\xi_i'')c_{hg,g^{-1}}c_{h,h^{-1}}c_e\xi_j').
\end{equation}

By (\ref{33}) and (\ref{44}), it suffices to show that
$$
c_{h,gh^{-1}g^{-1}}^{-1}h(\xi_i')h(\xi_j'')c_{hg^{-1},g}^{-1}hg^{-1}(\xi_i'')c_{hg^{-1},gh^{-1}g^{-1}}c_{g^{-1},g}c_e\xi_j'
$$
\begin{equation}\label{55}
=c_{hgh^{-1},g^{-1}}^{-1}hgh^{-1}(\xi_i')hgh^{-1}(\xi_j'')c_{hg,h^{-1}}^{-1}hg(\xi_i'')c_{hg,g^{-1}}c_{h,h^{-1}}c_e\xi_j'.
\end{equation}
Let us apply $h^{-1}$ to both hand sides of the latter equality. We will start with the left-hand side:
$$
\underline{h^{-1}(c_{h,gh^{-1}g^{-1}}^{-1})c_{h^{-1},h}c_e}\xi_i'\xi_j''c_e^{-1}\underline{c_{h^{-1},h}^{-1}h^{-1}(c_{hg^{-1},g}^{-1})c_{h^{-1},hg^{-1}}}g^{-1}(\xi_i'')\times
$$
$$
\times \underline{c_{h^{-1},hg^{-1}}^{-1}h^{-1}(c_{hg^{-1},gh^{-1}g^{-1}})}h^{-1}(c_{g^{-1},g}c_e\xi_j')
$$
$$
\stackrel{by (\ref{cocycle})}=c_{h^{-1},hgh^{-1}g^{-1}}\xi_i'\underline{\xi_j''c_e^{-1}c_{g^{-1},g}^{-1}}g^{-1}(\xi_i'')c_{g^{-1},gh^{-1}g^{-1}}c_{h^{-1},g^{-1}}^{-1}h^{-1}(\underline{c_{g^{-1},g}c_e\xi_j'})
$$
$$
\stackrel{by (\ref{inter})}=c_{h^{-1},hgh^{-1}g^{-1}}\xi_i'\xi_j''g^{-1}(\xi_i'')\underline{c_{g^{-1},gh^{-1}g^{-1}}c_{h^{-1},g^{-1}}^{-1}}h^{-1}(\xi_j')
$$
$$
=c_{h^{-1},hgh^{-1}g^{-1}}\xi_i'\xi_j''g^{-1}(\xi_i'')g^{-1}(c_e^{-1}c_{g,g^{-1}}^{-1}g(c_{g^{-1},gh^{-1}g^{-1}}c_{h^{-1},g^{-1}}^{-1})c_{g,g^{-1}}c_e)h^{-1}(\xi_j')
$$
$$
=c_{h^{-1},hgh^{-1}g^{-1}}\underline{\xi_i'}\xi_j''g^{-1}(\underline{\xi_i''c_e^{-1}c_{g,g^{-1}}^{-1}g(c_{g^{-1},gh^{-1}g^{-1}}c_{h^{-1},g^{-1}}^{-1})c_{g,g^{-1}}c_e})h^{-1}(\xi_j')
$$
\begin{equation}\label{77}
\stackrel{by (\ref{inter})}=c_{h^{-1},hgh^{-1}g^{-1}}c_e^{-1}c_{g,g^{-1}}^{-1}g(c_{g^{-1},gh^{-1}g^{-1}}c_{h^{-1},g^{-1}}^{-1})c_{g,g^{-1}}c_e\xi_i'\xi_j''g^{-1}(\xi_i'')h^{-1}(\xi_j').
\end{equation}
Now we will apply $h^{-1}$ to the right-hand side of (\ref{55}):
$$
\underline{h^{-1}(c_{hgh^{-1},g^{-1}}^{-1})c_{h^{-1},hgh^{-1}}}gh^{-1}(\xi_i')gh^{-1}(\xi_j'')\underline{c_{h^{-1},hgh^{-1}}^{-1}h^{-1}(c_{hg,h^{-1}}^{-1})c_{h^{-1},hg}}g(\xi_i'')\times
$$
$$
\times
\underline{c_{h^{-1},hg}^{-1}h^{-1}(c_{hg,g^{-1}})h^{-1}(c_{h,h^{-1}})h^{-1}(c_e)}h^{-1}(\xi_j')
$$
$$
\stackrel{by(\ref{cocycle})}=c_{h^{-1},hgh^{-1}g^{-1}}c_{gh^{-1},g^{-1}}^{-1}gh^{-1}(\xi_i')\underline{gh^{-1}(\xi_j'')}c_{g,h^{-1}}^{-1}g(\xi_i'')c_{g,g^{-1}}c_eh^{-1}(\xi_j')
$$
$$
=c_{h^{-1},hgh^{-1}g^{-1}}c_{gh^{-1},g^{-1}}^{-1}gh^{-1}(\xi_i')c_{g,h^{-1}}^{-1}g(\underline{h^{-1}(\xi_j'')})g(\xi_i'')c_{g,g^{-1}}c_e\underline{h^{-1}(\xi_j')}
$$
$$
\stackrel{by(\ref{Ginv})}=c_{h^{-1},hgh^{-1}g^{-1}}c_{gh^{-1},g^{-1}}^{-1}gh^{-1}(\xi_i')c_{g,h^{-1}}^{-1}g(\underline{c_{h,h^{-1}}c_e\xi_j''c_e^{-1}c_{h,h^{-1}}^{-1}})g(\xi_i'')
c_{g,g^{-1}}c_e\underline{c_{h,h^{-1}}c_e\xi_j'c_e^{-1}c_{h,h^{-1}}^{-1}}
$$
$$
\stackrel{by(\ref{inter})}=c_{h^{-1},hgh^{-1}g^{-1}}c_{gh^{-1},g^{-1}}^{-1}\underline{gh^{-1}(\xi_i')c_{g,h^{-1}}^{-1}}g(\xi_j'')g(\xi_i'')c_{g,g^{-1}}c_e\xi_j'
$$
$$
=c_{h^{-1},hgh^{-1}g^{-1}}c_{gh^{-1},g^{-1}}^{-1}c_{g,h^{-1}}^{-1}\underline{g(h^{-1}(\xi_i'))g(\xi_j'')}g(\xi_i'')c_{g,g^{-1}}c_e\xi_j'
$$
$$
=c_{h^{-1},hgh^{-1}g^{-1}}c_{gh^{-1},g^{-1}}^{-1}c_{g,h^{-1}}^{-1}g(\underline{h^{-1}(\xi_i')\xi_j''})g(\xi_i'')c_{g,g^{-1}}c_e\underline{\xi_j'}
$$
$$
\stackrel{by (\ref{inter})}=c_{h^{-1},hgh^{-1}g^{-1}}c_{gh^{-1},g^{-1}}^{-1}c_{g,h^{-1}}^{-1}\underline{g(\xi_j'')g(\xi_i''})c_{g,g^{-1}}c_e\xi_j'h^{-1}(\xi_i')
$$
$$
=c_{h^{-1},hgh^{-1}g^{-1}}c_{gh^{-1},g^{-1}}^{-1}c_{g,h^{-1}}^{-1}g(\xi_j''\xi_i'')c_{g,g^{-1}}c_e\xi_j'h^{-1}(\xi_i')
$$
$$
\stackrel{by (\ref{inter})}=c_{h^{-1},hgh^{-1}g^{-1}}c_{gh^{-1},g^{-1}}^{-1}c_{g,h^{-1}}^{-1}\underline{g(\xi_j'')}c_{g,g^{-1}}c_e\xi_i''\underline{\xi_j'}h^{-1}(\xi_i')
$$
$$
\stackrel{by (\ref{Ginv})}=c_{h^{-1},hgh^{-1}g^{-1}}c_{gh^{-1},g^{-1}}^{-1}c_{g,h^{-1}}^{-1}c_{g,g^{-1}}c_e\xi_j''\xi_i''g^{-1}(\xi_j')h^{-1}(\xi_i')
$$

Now let us compare the latter expression with (\ref{77}): clearly, to finish proving the torus axiom, it remains to show that
$$
c_e^{-1}c_{g,g^{-1}}^{-1}g(c_{g^{-1},gh^{-1}g^{-1}}c_{h^{-1},g^{-1}}^{-1})
=c_{gh^{-1},g^{-1}}^{-1}c_{g,h^{-1}}^{-1}.
$$
By (\ref{cocycle}), the right-hand side equals $c_{g,h^{-1}g^{-1}}^{-1}g(c_{h^{-1},g^{-1}}^{-1})$, so the equality reduces to
$$
c_e^{-1}c_{g,g^{-1}}^{-1}g(c_{g^{-1},gh^{-1}g^{-1}})
=c_{g,h^{-1}g^{-1}}^{-1}
$$
or, equivalently,
$$
g(c_{g^{-1},gh^{-1}g^{-1}})c_{g,h^{-1}g^{-1}}
=c_{g,g^{-1}}c_e.
$$
The latter equality follows from (\ref{cocycle}). The torus axiom is proved.

\bigskip

\small{\tt Mathematics Department, Kansas State University}

\small{\tt e-mail: shklyarov@math.ksu.edu}
\end{document}